\documentclass[a4paper,12pt]{article}

\usepackage[T1]{fontenc}
\usepackage{amsmath, graphicx, color, amsthm, centernot, amssymb}
\usepackage[font=small,labelfont=bf,width=12.5cm]{caption}

\newtheorem{claim}{Claim}
\newtheorem{theorem}{Theorem}

\newtheorem{conjecture}{Conjecture}

\newtheorem{problem}{Problem}
\newtheorem{question}{Question}
\newtheorem{observation}{Observation}

\newcommand{\G}{\ensuremath{\vec{G}}}
\newcommand{\inter}{\ensuremath{\mathrm{int}}}
\newcommand{\Inter}{\ensuremath{\mathrm{Int}}}
\newcommand{\exter}{\ensuremath{\mathrm{ext}}}

\newenvironment{proofclaim}[1][]%
    {\noindent \emph{Proof.} {}{#1}{}}{$~$\hfill $~\blacklozenge$ \vspace{0.2cm}}

\title{Note on $3$-Choosability of Planar Graphs with Maximum Degree $4$}
\author
{
	Fran\c{c}ois Dross\thanks{LIRMM, Universit\'{e} de Montpellier \& CNRS, France. 
		E-Mail: \texttt{dross@lirmm.fr}}, \
	Borut Lu\v{z}ar\thanks{Faculty of Information Studies, Novo mesto, Slovenia.
		E-Mail: \texttt{borut.luzar@gmail.com}} \footnotemark[3], \
	M\'{a}ria Macekov\'{a}\thanks{Faculty of Science, Pavol Jozef \v Saf\'{a}rik University, Ko\v{s}ice, Slovakia.
		E-Mail: \texttt{\{maria.macekova,roman.sotak\}@upjs.sk}}, \
	Roman Sot\'{a}k\footnotemark[3]
}

\begin{document}
\maketitle

{
\abstract
{
	Deciding whether a planar graph (even of maximum degree $4$) is $3$-colorable is NP-complete.
	Determining subclasses of planar graphs being $3$-colorable has a long history, but since Gr\"{o}tzsch's result
	that triangle-free planar graphs are such, most of the ef\mbox{}fort was focused to solving Havel's and Steinberg's conjectures.
	In this paper, we prove that every planar graph obtained as a subgraph of the medial graph of any bipartite plane graph is $3$-choosable.	
	These graphs are allowed to have close triangles (even incident), and have no short cycles forbidden, 
	hence representing an entirely different class than the graphs inferred by the above mentioned conjectures. 
}

\bigskip
{\noindent\small \textbf{Keywords:} medial graph, plane graph, $3$-colorability, $3$-choosability, Alon-Tarsi Theorem.}
}

\section{Introduction}

In this paper we consider the problem of $3$-colorability of a subclass of $4$-regular planar graphs.
The initial motivation came from the following problem about medial graphs (defined in Section~\ref{sec:prel}), 
proposed by Czap, Jendro\v{l}, and Voigt~\cite[Problem~3.9]{CzaJenVoi18}.
\begin{problem}
	Is there a bipartite plane graph $G$ such that its medial graph has chromatic number $4$?
\end{problem}
\noindent  
We prove that such a graph does not exist. Even more, the medial graphs of bipartite plane graphs are $3$-choosable.
Note that one cannot omit the bipartiteness condition as already the medial graph of the complete graph on four vertices 
with one edge subdivided is not $3$-colorable.

Medial graphs are planar and $4$-regular, 
hence the problem reduces to investigating $3$-colorability of a subclass of planar graphs with maximum degree $4$.
While every connected graph of maximum degree $3$ which is not isomorphic to $K_4$ is $3$-colorable by Brooks' Theorem, 
already in the class of planar graphs with maximum degree $4$ deciding whether a graph admits a $3$-coloring
is NP-complete~\cite{GarJohSto76}.
Due to this fact, and even more due to the famous Four Color Theorem, 
the problem of $3$-coloring received a lot of attention in the class of planar graphs. 
For any plane triangulations, Heawood found a necessary and sufficient condition~\cite{Hea98} showing that
it is $3$-colorable if and only if all its vertices have even degrees.
Generalizations of this statement have been given in~\cite{DikKowKur02,EllFleKocWen04} and just recently in~\cite{Koc18}.

On the other hand, a well-known result by Gr\"{o}tzsch's~\cite{Gro59} shows that if there are no cycles of length $3$ in a planar graph, then it is $3$-colorable. 
This result was later improved by Gr\"{u}nbaum~\cite{Gru63} to planar graphs with at most three triangles.
The original proof was faulty, but was later corrected by Aksenov~\cite{Aks74}. 
Finally, a simpler proof of his result was recently given in~\cite{BorKosLidYan14}.

Allowing some triangles in a graph, but still retaining $3$-colorability yielded two intriguing conjectures.
First, Havel~\cite{Hav69} conjectured that a $3$-colorable planar graph may contain many triangles as long
as they are sufficiently far apart. 
This conjecture was recently proved by Dvo\v{r}\'{a}k, Kr\'{a}\v{l}, and Thomas~\cite{DvoKraTom16}
(they announced the result already in 2009, but included it in a series of papers on $3$-colorability 
of triangle-free graphs on surfaces; cf.~\cite{DvoKraTom16b}).

The second conjecture is due to Steinberg~\cite{Ste93}. It allows arbitrarily many triangles but it forbids short cycles.
Namely, Steinberg conjectured that every planar graph without cycles of length $4$ and $5$ is $3$-colorable.
The conjecture was disproved by Cohen-Addad et al.~\cite{CohHebKraLiSal17};
however a number of weaker results have been proved, perhaps the closest
being due to Borodin et al.~\cite{BorGleMonRas09}, stating
that every planar graph without cycles of length $5$ and $7$, and without adjacent triangles is $3$-colorable 
(see also~\cite{BorGleRasSal05,BorKosLidYan14,BorMonRas10} for other results on this conjecture).

All the problems listed above are even harder in a more general setting of list coloring.
As shown by Voigt~\cite{Voi93}, planar graphs are not $4$-choosable; 
Thomassen~\cite{Tho94} found a beautiful proof that they are $5$-choosable. 
It is a rather simple observation that planar bipartite graphs are not $2$-choosable,
and a bit harder that they are $3$-choosable~\cite{AloTar92}.
Therefore it is not surprising that an equivalent of Gr\"{o}tzsch's result does not hold in this setting; 
as shown by Voigt~\cite{Voi95}, there are triangle-free planar graphs which are not $3$-choosable.
Thomassen~\cite{Tho95} however proved that girth $5$ is a sufficient condition for their $3$-choosability. 

An analogue of Havel's conjecture hence requires cycles of length $3$ and also $4$ to be sufficiently distant.
Dvo\v{r}\'{a}k~\cite{Dvo14} proved that distance $26$ between them is sufficient, 
while the best lower bound requires distance at least $4$~\cite{AksMel80}.
The list-version of Steinberg's conjecture clearly requires more excluded cycles. 
So far, there has been a number of partial results towards solving it (cf.~\cite{DvoLidSkr09} and references therein for more details), 
with currently the best one, due to Dvo\v{r}\'{a}k and Postle~\cite{DvoPos18}, showing that planar graphs without cycles of lengths from $4$ to $8$
are $3$-choosable. It is still not known if it suffices to forbid only cycles from $4$ to $7$ or even from $4$ to $6$.

We contribute to the above described rich field of research with the following theorem, which does assume a special structure of a graph,
but does not particularly bound the number of triangles. 
The triangles can even have common vertices, so the distance between them can be arbitrarily small.
On the other hand, the graph can contain cycles of any length.
\begin{theorem}
	\label{thm:main}
	Every medial graph (with eventual loops removed) of a bipartite plane graph is $3$-choosable.
\end{theorem}

We prove the theorem in Section~\ref{sec:main}, while in Section~\ref{sec:prel}, 
we introduce notation and auxiliary results. We conclude with some open problems in Section~\ref{sec:conc}.

\section{Preliminaries}
\label{sec:prel}

For a plane graph $G$, with $V(G)$, $E(G)$, and $F(G)$ we denote its set of vertices, edges, and faces, respectively.
The \textit{medial graph} $M(G)$ of $G$ is the graph with the vertex set $V(M(G)) = E(G)$, two vertices $u$ and $v$
being adjacent if the edges of $G$ corresponding to $u$ and $v$ appear successively on the boundary of some face of $G$.
Thus, every medial graph is $4$-regular.

By the \textit{intersection} of graphs $G_1$ and $G_2$, denoted $G_1 \cap G_2$, we mean the graph with 
the sets of vertices and edges that are contained in both graphs. 

\begin{figure}[ht!]
	$$
		\includegraphics{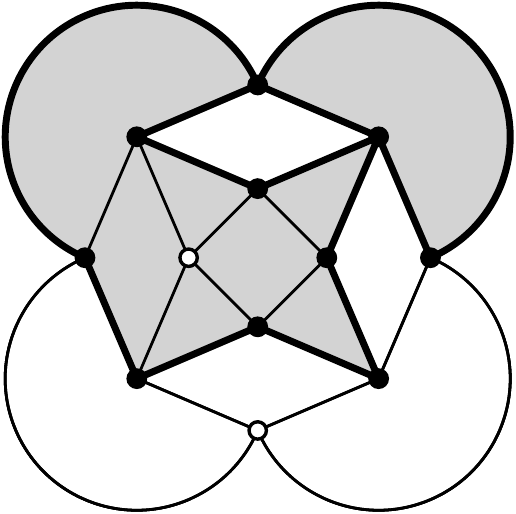}
	$$
	\caption{A graph $G$ (all the vertices and edges) and its Eulerian subgraph $H$ represented with full vertices and thick edges.
		The graph $H$ has four faces; the two non-shaded represent green faces, and the two shaded represent blue faces.}
	\label{fig:ex-inout}
\end{figure}

\begin{figure}[ht!]
	$$
		\includegraphics{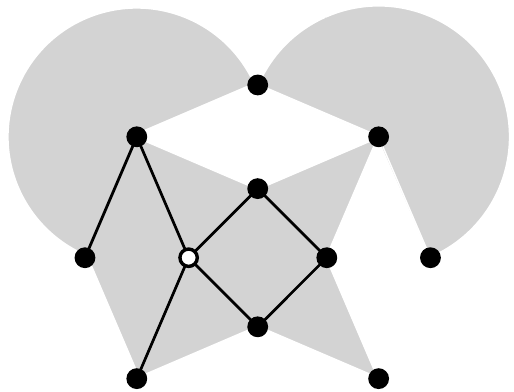} \quad \quad
		\includegraphics{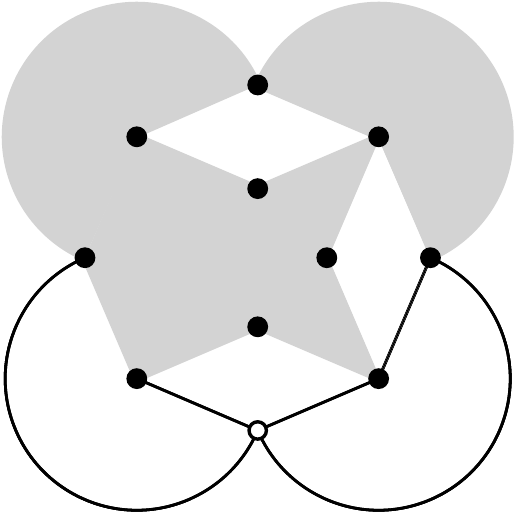}
	$$
	\caption{On the left, $\inter(H)$ is shown (we retain the initial colors of faces, although some of their boundary edges are removed).
		Similarly, on the right, the graph $\exter(H)$ is depicted.}
	\label{fig:ex-inout2}
\end{figure}

Given a plane graph $G$,
we define the boundary, interior, and exterior of any plane Eulerian subgraph $H$ of $G$ in the following way.
First, color the faces of $H$ with two colors such that adjacent faces receive distinct colors 
(this is possible as the dual of $H$ is bipartite).
Let the outer face of $H$ be colored green and its adjacent faces blue
(see Figure~\ref{fig:ex-inout} for an example).
The boundary of $H$, $\partial(H)$, is the graph $H$ itself.
The interior $\inter(H)$ is the graph induced by the vertices of $G$ lying in the blue faces of $H$
together with the vertices of $H$ without the edges of $H$,
and 
the exterior $\exter(H)$ is the graph induced by the vertices of $G$ lying in the green faces of $H$
together with the vertices of $H$ without the edges of $H$ (see Figure~\ref{fig:ex-inout2}).

Similarly, for a subgraph $X$ of $G$, we define
$$
	\partial_X(H) = \partial(H) \cap X, \quad \inter_X(H) = \inter(H) \cap X, \quad \exter_X(H) = \exter(H) \cap X.
$$

We say that $L$ is a \textit{list-assignment} for the graph $G$ if it assigns a list $L(v)$ of possible colors to each vertex $v$ of $G$. 
If $G$ admits a proper coloring $\varphi_l$ such that $\varphi_l(v) \in L(v)$ for all vertices in $V(G)$, 
then we say that $G$ is \textit{$L$-colorable} or $\varphi_l$ is an {\em $L$-coloring} of $G$. 
The graph $G$ is \textit{$k$-choosable} if it is $L$-colorable for every assignment $L$, where $|L(v)| \ge k$ for every $v \in V(G)$. 
The \textit{list chromatic number} $\chi_{l}(G)$ of $G$ is the smallest $k$ such that $G$ is $k$-choosable.

For a directed graph (digraph) $D$, we define the \textit{indegree} (denoted $d^-(v)$) and \textit{outdegree} (denoted $d^+(v)$) of a vertex $v$
as the number of directed edges having $v$ as a terminal and an initial vertex, respectively.
A subdigraph $H$ of a digraph $D$ is called 
\textit{Eulerian} if the indegree $d^{-}_{H}(v)$ of every vertex $v$ of $H$ is equal to its outdegree $d^{+}_{H}(v)$. 
The digraph $H$ is \textit{even} if it has even number of directed edges, otherwise, it is \textit{odd}. 
Let $E^e(D)$ and $E^o(D)$ be the numbers of even and odd spanning Eulerian subgraphs of $D$, respectively.

The following well-known result due to Alon and Tarsi~\cite{AloTar92} plays the key role in our proof of Theorem~\ref{thm:main}.
\begin{theorem}[Alon \& Tarsi, 1992]
	\label{thm:AT}
	Let $D$ be a directed graph, and let $L$ be a list-assignment such that 
	$|L(v)| \ge d^{+}_{D}(v) + 1$ for each $v \in V(D)$. 
	If $E^e(D) \neq E^o(D)$, then $D$ is $L$-colorable. 
\end{theorem}

\section{Proof of Theorem~\ref{thm:main}}
\label{sec:main}

In this section, we present a proof of Theorem~\ref{thm:main} using the result of Alon and Tarsi mentioned in the previous section.
We first discuss the assumptions of Theorem~\ref{thm:main}. The medial graph of a graph with a vertex $v$ of degree $1$ contains a loop,
since the boundary of the face containing $v$ uses the edge incident to $v$ twice. As a graph with loops cannot admit a proper vertex coloring, 
since a vertex incident with a loop is adjacent to itself, we rather consider its subgraph with the loops removed.

\begin{proof}
	Let $B$ be a bipartite plane graph. 
	For technical reasons described above, we recursively remove all its vertices of degree at most $1$ to obtain the graph $B'$
	with minimum degree at least $2$. 	
	Let $G$ be the medial graph of $B'$ (notice that $G$ may have parallel edges if there are vertices of degree $2$ in $B'$).
	Let $M'(B)$ denote $M(B)$ with loops removed.
	Fix a list-assignment $L$ such that $|L(v)| = 3$ for all $v \in V(M'(B))$.	 
	Below we describe a procedure to $L$-color $G$. 
	It is straightforward to observe that any $L$-coloring of $G$ can be extended to an $L$-coloring of $M'(B)$, 
	since every vertex to be colored has at most two colored neighbors.

	There are two types of faces in $G$: the ones corresponding to the vertices of $B$
	(we call them \textit{black}), and the ones corresponding to the faces of $B$ 
	(we call them \textit{white}). 
	Notice that all white faces have even length, since $B$ is bipartite. 
	Moreover, every edge in $G$ is incident to two faces, one black and one white (see the left graph in Figure~\ref{fig:ex1} for an example).
	\begin{figure}[ht!]
		$$
			\includegraphics{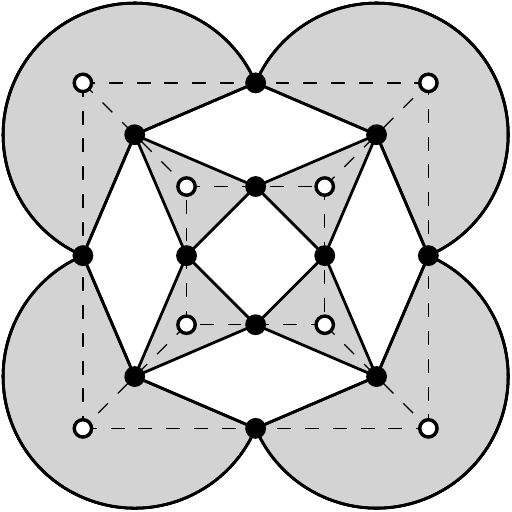} \quad \quad
			\includegraphics{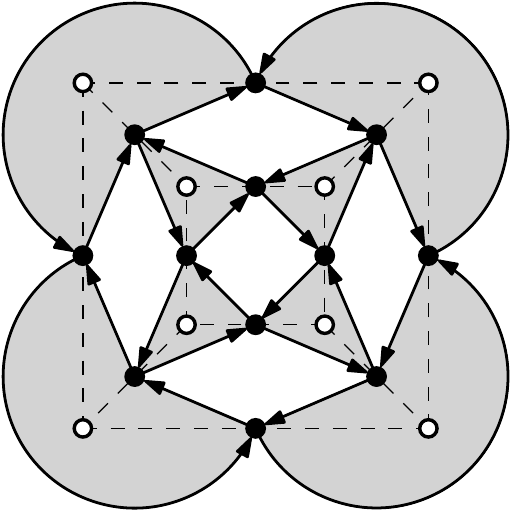}
		$$
		\caption{The left graph is the medial graph $M(Q_3)$ of the cube $Q_3$; its vertices and edges are depicted as full circles and solid edges,
			while the vertices and edges of $Q_3$ are depicted with empty circles and dashed edges. The black faces of $M(Q_3)$ are shaded.
			The right graph is the directed graph $\protect\overrightarrow{M(Q_3)}$ such that the edges have a black face always on their left hand side.}
		\label{fig:ex1}
	\end{figure}
	
	Let $\G$ be a directed graph obtained from $G$ by directing the edges 
	such that each edge has its black face on the left hand side when going from 
	its initial vertex to its terminal vertex. 
	This in particular means that every vertex has precisely two incoming and two outgoing edges,
	and therefore $d^-(v) = d^+(v) = 2$ for every $v \in V(\G)$.
	Apart from the regularity, we will also use the following fact,
	implied by the choice of orientation.
	\begin{observation}
		\label{obs:cut1}
		Two consecutive edges on a directed cycle always appear successively on the boundary of some face.
	\end{observation}
	In other words, the choice of orientation guarantees that an incoming edge at some vertex is facially adjacent only 
	to outgoing edges.
	
	By Observation~\ref{obs:cut1}, we also infer the following relationship between two directed cycles.
	\begin{observation}
		\label{obs:cut}
		Let $D_1$ and $D_2$ be two directed cycles in $\G$ intersecting (i.e., having some common vertices) 
		in such a way that
		$E(\partial(D_2) \cap \inter(D_1)) \neq \emptyset$ and $E(\partial(D_2) \cap \exter(D_1)) \neq \emptyset$.
		Then $E(D_1) \cap E(D_2) \neq \emptyset$.
	\end{observation}
	Namely, if $D_2$ has edges in the interior and in the exterior of $D_1$, then they also share edges.
				
	In what follows, we will show that every odd Eulerian spanning subgraph of $\G$ 
	can be injectively mapped to an even Eulerian spanning subgraph of $\G$. 
	We will also show there is an even Eulerian subgraph of $\G$ to which no
	odd subgraph is mapped, and thus fulf\mbox{}ill the assumptions of Theorem~\ref{thm:AT}.
	That will imply $3$-choosability of $G$.
	
	We distinguish two types of directed cycles in $\G$: 
	by the definition and Observation~\ref{obs:cut1}, all the edges of a given directed cycle $C$
	are incident either to black faces or to white faces in the interior of $C$.
	We refer to the former as \textit{black cycles} (see Figure~\ref{fig:excyc} for an example) and to the latter as \textit{white cycles}.	
	\begin{figure}[ht!]
		$$
			\includegraphics{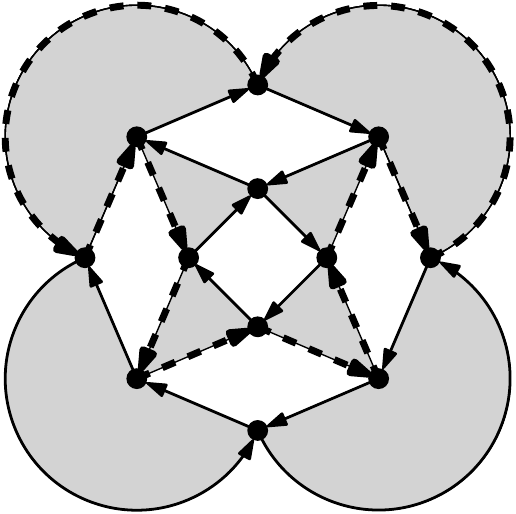}
		$$
		\caption{A black cycle in $\protect\overrightarrow{M(Q_3)}$ with its edges depicted by fat dashed lines.}
		\label{fig:excyc}
	\end{figure}	
	Similarly, we say that an Eulerian graph is \textit{white} if it is comprised of white cycles only.
	
	For a graph $X$, we denote its complement by $\overline{X}$.
	For a cycle $D$, the \textit{$D$-complement} of a spanning Eulerian subgraph $X$ of $\G$
	is the spanning Eulerian subgraph $\overline{X}^D$ with the edge set 
	$$
		E(\overline{X}^D) = E(\exter_X(D)) \cup E(\inter_{\overline{X}}(D)) \cup E(\partial_{\overline{X}}(D))\,.
	$$	
	The fact that $\overline{X}^D$ is also Eulerian follows from Observation~\ref{obs:cut1}.

	\begin{claim}
		\label{cl:evenodd}
		For an odd black cycle $D$, the $D$-complement of an odd (even) Eulerian spanning subgraph $X$ is an even (odd) Eulerian spanning subgraph $\overline{X}^D$.
	\end{claim}
	
	\begin{proofclaim}
		All the edges of $\inter(D)$ comprise an edge-disjoint union of cycles around white faces 
		(since $D$ is a black cycle it means that all the white faces in $\inter(D)$ have all boundary edges in $\inter(D)$ also).
		Since every white face has even length, there is an even number of edges in $\inter(D)$.
		Hence, the parity of the number of edges in $\inter_X(D)$ is the same as the parity of the number of edges in $\inter_{\overline{X}}(D)$. 
		So the parity of the number of edges of $X$ is different from the parity of the number of edges of its $D$-complement
		as one of them contains an odd number of edges of $D$, and the other one an even number of edges of $D$.
		(Recall that $E(\exter(\overline{X}^D)) = E(\exter(X))$.)
	\end{proofclaim}
	
	Before we proceed, we prove a claim about cycles in the intersection of an Eulerian graph and its $D$-complement.	
	\begin{claim}
		\label{cl:int}
		Let $X$ be an Eulerian spanning subgraph of $\vec{G}$, 
		and let $D$ be an odd white Eulerian subgraph of $X$. 
		Then, there is an odd black cycle in $\inter_X(D)$ or $\inter_{\overline{X}^D}(D)$.
	\end{claim}

	\begin{proofclaim}
		Suppose the contrary and let $D$ be minimal in terms of the number of faces in its interior.
		On the internal side of the edges of the boundary $\partial(D)$, there are white faces (by definition).		
		Take the edges of these faces which are not in $\partial(D)$.
		There is also an odd number of them, since white faces have even length.
		At least one such edge $e$ is in $X$, otherwise there is an odd black cycle in $\inter_{\overline{X}^D}(D)$.
		But, $e$ also belongs to some black cycle $C$ in $\inter_X(D)$. If $C$ is even, we add it to $D$ (in which case it is white from the point of view of $D$),
		removing the edges of $\inter(C)$ from $D$ and obtaining a smaller graph, which contradicts the minimality of $D$.
		Otherwise, $C$ is an odd black cycle, and hence the claim is established.
	\end{proofclaim}
	
	Denote by $\mathcal{E}$ the set of all Eulerian spanning subgraphs of $\G$.
	Let $\mathcal{O}$ be a sorted set of all odd black cycles in $\G$,
	sorted in ascending order by the number of faces they contain in their interiors.
	Suppose there are $k$ cycles, $C_1,C_2,\dots,C_k$, in $\mathcal{O}$.
	For every $i$, $1 \le i \le k$, in consecutive order, 
	\textit{we repeatedly remove all $X \in \mathcal{E}$ which either contain all the edges of $C_i$ or none of them}.
	Hence, if in the step $i$ we remove from $\mathcal{E}$ some $X$, then we also remove its $C_i$-complement if it is still 
	in $\mathcal{E}$. In the proof of the claim below, we show that such pairs are always removed at the same step.
	
	\begin{claim}
		The number of odd Eulerian spanning subgraphs removed from $\mathcal{E}$ at step $i$
		is equal to the number of even such subgraphs.
	\end{claim}
	
	\begin{proofclaim}
		By Claim~\ref{cl:evenodd}, we have that an Eulerian spanning subgraph $X$ is of different parity from
		 its $C_i$-complement.
		Hence, proving that removing $X$ from $\mathcal{E}$ at step $i$ implies removal of its $C_i$-complement at step $i$
		establishes the claim.
		
		Suppose the contrary, and let $i$ be minimal such that there is some $X$ in $\mathcal{E}$ 
		whose $C_i$-complement $\overline{X}^{C_i}$ is not in $\mathcal{E}$, 
		i.e. has been removed in some step $j$, with $j<i$. Then $\overline{X}^{C_i}$ contains all the edges of $C_j$ or none of them.
		First, notice that $C_j$ is not completely contained in $\partial(C_i) \cup \inter(C_i)$,
		since, by the definition of $\overline{X}^{C_i}$, $X$ would also contain all the edges of $C_j$ or none of them,
		meaning that $X$ would also be removed at the step $j$.
		Hence, $C_j$ has an edge in $\exter(C_i)$. 
		For the same reason, $C_j$ is not completely contained in $\exter(C_i)$, and thus has an edge in $\partial(C_i) \cup \inter(C_i)$.
		In fact, by Observation~\ref{obs:cut}, it follows that there must be some edge of $C_j$ on $\partial(C_i)$.
		Finally, there is also some edge of $C_j$ in $\inter(C_i)$, otherwise $C_i$ is either contained in $\partial(C_j) \cup \inter(C_j)$,
		in which case $i < j$, or $C_j$ is white. In both cases, we obtain a contradiction.
		
		Now, we show that there is some odd black cycle in $\partial(C_i) \cup \inter(C_i)$ (distinct from $C_i$)
		such that $X$ either contains all its edges or none of them.
		Suppose first that $X$ contains all the edges of $C_i$.
		Since, by the above argumentation, $C_i$ and $C_j$ have some edges in common, 
		this implies that none of the edges of $C_j$ is in $\overline{X}^{C_i}$, and so all the edges of $C_j$ 
		in $\partial(C_i) \cup \inter(C_i)$ are also in $X$.		
		
		By Observation~\ref{obs:cut1}, $X_{\Inter} = \partial(C_i) \cup \inter_X(C_i)$ is Eulerian with the outer face of odd length.
		By the Handshake Lemma, there is also an innerface $f$ of odd length (note that $f$ is a face of $X_{\Inter}$, 
		but not necessarily of $\partial(C_i) \cup \inter(C_i)$). 
		Moreover, $f$ is not bounded by $C_i$, since, by the argumentation above, there is at least one edge in $\inter_X(C_i)$. 
		Let $C$ be the cycle bounding $f$ in $X_{\Inter}$. 
		If $C$ is black, we are done, since there is some $\ell < i$ such that $C = C_\ell$, and so $X$ would be removed at the step $\ell$.
		
		Otherwise, $C$ is white and by Claim~\ref{cl:int}, 
		there is some odd black cycle $C_\ell \in \mathcal{O}$ in $\inter(C)$, 
		such that $\overline{X}^{C_i}$ contains all the edges of $C_\ell$ or none of them, and $X$ contains none of the edges of $C_\ell$ or all of them, 
		respectively, which means that $X$ would be removed at the step $\ell$.
		
		Suppose now that $X$ contains none of the edges of $C_i$, and therefore $\overline{X}^{C_i}$ contains all of them. 
		In this case, we use analogous argumentation as in the previous paragraph that there is an odd black cycle such that 
		$\partial(C_i) \cup \inter_{\overline{X}^{C_i}}(C_i)$ contains either all or none of its edges, 
		meaning that $X$ contains none or all of its edges, respectively, which again
		implies that $X$ would have been removed in some of the previous steps, a contradiction.
	\end{proofclaim}
	
	In every odd Eulerian graph there is an odd cycle. Since in $\G$ every directed cycle is either black or white,
	by Claim~\ref{cl:int}, every odd Eulerian spanning subgraph of $\G$ contains an odd black cycle or a complement
	of an odd black cycle. This implies that after all cycles from $\mathcal{O}$ are removed, 
	there is no odd Eulerian spanning subgraph left in $\mathcal{E}$.	
	
	However, there is at least one even Eulerian spanning subgraph, which contains at least one edge of every odd black cycle in $\G$,
	but not all edges of any. We guarantee its existence by the following claim.
	\begin{claim}	
		White faces of $G$ can be colored with two colors, red and blue, such that every odd black cycle
		shares an edge with the boundary of at least one red and at least one blue face.
	\end{claim}
	\begin{proofclaim}
		Let $H$ be the graph whose vertex set is formed by the white faces of $G$ and two vertices of $H$ 
		are connected if the corresponding white faces share a vertex in $G$. The graph $H$ is planar, 
		and hence we can color its vertices with four colors, say 1, 2, 3, and 4, by the Four Color Theorem.
		Color the white faces of $G$ whose corresponding vertices in $H$ are colored with $1$ or $2$, with red,
		and the other white faces with blue. 
		
		Let $C$ be an odd black cycle. 
		By the orientation of the graph, every vertex of $C$ has its two incident edges 
		that are not in $C$ either both in $\inter(C)$ or both in $\exter(C)$. 
		Let $V_i$ be the set of the vertices of $C$ that have two incident edges in $\inter(C)$ 
		and let $V_e$ be the set of vertices of $C$ that have two incident edges in $\exter(C)$. 
		As every vertex has degree $4$, the number of edges in $\inter(C)$ is equal to twice the number of vertices 
		that have all of their edges in $\inter(C)$, plus $|V_i|$. 
		As the edges in $\inter(C)$ are the edges of the disjoint union of the boundaries of white faces, 
		and as white faces are even, 
		there are even number of edges in $\inter(C)$, thus $|V_i|$ is even. 
		The number of white faces that share an edge with $C$ is equal to $|V_e|$, and as $C$ is odd, $|V_e|$ is odd. 
		Therefore there is an odd number of white faces that share edges with $C$, 
		and these form an odd cycle in $H$. 
		As there are three colors needed for coloring an odd cycle, the claim is established.
	\end{proofclaim}
	
	By taking the edges of the union of the boundaries of the red faces, 
	we obtain an even Eulerian subgraph that contains at least one edge of every odd black cycle in $\G$, 
	but not all edges of any. 
	This even eulerian subgraph is still in $\mathcal{E}$.
	
	Hence, we have proved that there are more even Eulerian spanning subgraphs in $\G$ than
	odd Eulerian spanning subgraphs and thus fulfill the assumptions of Theorem~\ref{thm:AT}.
	This means that $G$ is $3$-choosable.
\end{proof}

\section{Discussion and Further Work}
\label{sec:conc}

In this paper, we answered the question of Czap, Jendro\v{l}, and Voigt~\cite[Problem~3.9]{CzaJenVoi18} 
about chromatic number of medial graphs of bipartite plane graphs. We used an application of the Theorem of Alon and Tarsi
and to satisfy the main assumption of it, we strongly used the fact our graphs have maximum degree $4$ and that its
faces can be properly colored in two colors, where one color class contains only even faces. 
It is not clear if the former condition is really needed.
In fact, we believe that the following conjecture can be answered in affirmative.
\begin{conjecture}
	\label{conj:main}
	Every simple plane graph whose faces can be properly colored with two colors 
	such that one color class contains only even faces is $3$-colorable.	
\end{conjecture}
In Conjecture~\ref{conj:main}, we require a simple graph, since any plane graph with every edge replaced by two parallel edges
satisfies the assumption of being $2$-face-colorable, with the faces created by the parallel edges being even, and so any plane graph 
with chromatic number $4$ would be a counterexample. Our result, however, allows parallel edges since the maximum degree is limited to $4$.

Conjecture~\ref{conj:main} can be strengthened to the choosability version, but we are less certain about the answer.
\begin{question}
	Is every simple plane graph whose faces can be properly colored with two colors 
	such that one color class contains only even faces also $3$-choosable?	
\end{question}

Regarding the above conjecture and question, we would like to mention another somehow 
related direction of research which was pointed to us by a referee.
It is well known that the Four Color Theorem is equivalent to the fact that 
every bridgeless cubic planar graph has chromatic index $3$.
A little less known is that bridgeless cubic planar graphs are $3$-edge-choosable.
This was proved by Ellingham and Goddyn~\cite{EllGod96} 
(they also mention that this fact was reported by Jaeger and Tarsi in personal communication) 
using the Theorem of Alon and Tarsi.
In fact they proved even more, every $d$-regular planar graph of class I is $d$-edge-choosable.
Here we focus to the result on cubic graphs, since the medial graph of a cubic plane graph $G$
is just the line graph of $G$, and hence the result of Ellingham and Goddyn is somehow
of the same flavor as ours. 
The result does not assume that the original graph $G$ is bipartite, but it assumes bounded maximum degree.
The medial graph of $G$ again can have the faces properly colored with two colors, and the faces of one color (this time the black ones)
are all of length $3$ (in our case all white faces have even lengths).

This leads to a further question: if after coloring the faces of the medial graph with $2$ colors all faces of one color 
have lengths of the same parity, is it true that the medial graph is $3$-colorable ($3$-choosable)~\cite{Kai18}?
The question has a negative answer, consider e.g. the wheel on five spokes $W_6$. 
Its medial graph has all black faces of odd lengths, but it is not $3$-colorable (we leave it to the reader to verify this).
It does, however, remain the following.
\begin{question}
	Let $G$ be a simple plane graph whose faces can be properly colored with two colors
	such that one color class contains faces with lengths of the same parity.
	Is it true that from $3$-colorability of $G$ it follows that $G$ is also $3$-choosable?
\end{question}

\begin{figure}[htp!]
	$$
		\includegraphics{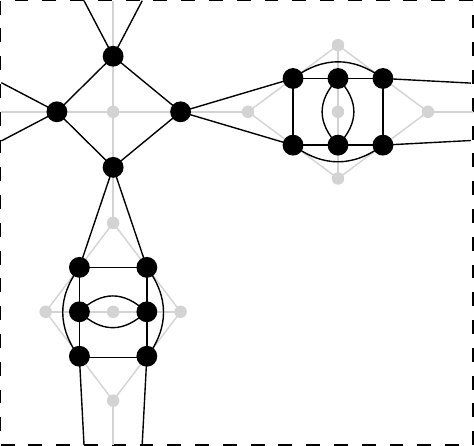}
	$$
	\caption{A medial graph (depicted in black) of a bipartite graph embedded on torus (depicted in grey) having chromatic number $4$.}
	\label{fig:ex-torus}
\end{figure}

Finally, let us discuss the planarity condition as one might ask, why limit only to plane graphs. 
Theorem~\ref{thm:main} does not hold in general for graphs that embed to other surfaces.
In Figure~\ref{fig:ex-torus}, we present a graph which needs $4$ colors for a vertex-coloring 
(as it is straightforward to check that $3$ colors do not suffice to color it, we leave it to the reader).

\paragraph{Acknowledgment.} 

The authors are indebted to a reviewer who pointed out the result about $3$-edge-choosability of cubic planar graphs. 
The first author was partially supported by the ANR grant HOSIGRA (contract number ANR-17-CE40-0022-03).
The second author was partly supported by the Slovenian Research Agency Program P1--0383 and joint Austrian-Slovenian project N1-0057.
The latter two authors were supported by the Slovak Research and Development Agency under the Contract No. APVV--15--0116 and by the Science Grant Agency - project VEGA 1/0368/16.

\bibliographystyle{plain}
{\small
	\bibliography{MainBase}
}

\end{document}